\documentclass[12pt,twoside]{article}
\usepackage[all]{}
\usepackage {amssymb,latexsym,amsthm,amscd,amsmath}
\newtheorem*{theorem*}{Theorem}
\theoremstyle{definition}

\numberwithin{equation}{subsection}

\newcommand{\ignore}[1]{}

\newcommand{\mynote}[1]{}
\begin{document}
\setcounter{section}{0}
\title{\bf The cyclicity problem for Albert algebras}
\author{Maneesh Thakur\\ \small{Indian Statistical Institute, Stat-Math-unit}\\ \small{8th Mile Mysore Road, RV College Post}\\\small{Bangalore 560059, India}\\
  \small{e-mail: maneesh.thakur@gmail.com}
}
\date{}
\maketitle
\begin{abstract}
  \noindent
      {Addressing a long-standing question of Albert dating back to 1965, we show that every Albert division algebra has an isotope that contains a cyclic cubic subfield. }
  
\end{abstract}

\section{Introduction}
More than fifty years ago, Abraham Adrian Albert raised a question that in modern terminology may be expressed as follows.
\vskip1mm
\emph{Does every Albert division algebra (of characteristic not $2$ or $3$) contain a cyclic cubic subfield?}
\vskip1mm

This was perhaps motivated by a theorem of Wedderburn, which states that every degree $3$ central division algebra over a field contains a cyclic cubic subfield. Albert algebras are Jordan algebras of degree $3$ (see \cite{J}, Chap. IX). 
The aforementioned question, which derives its importance not only from Albert's own work in (\cite{A-1}, \cite{A-2}) on the subject, but also from Springer's theory of cyclic compositions (see \cite{KMRT}, \cite{P3}, \cite{SV0}, \cite{SV}), has been answered affirmatively by Petersson-Racine (\cite{PR1}, Thm. 4) in the presence of the cube roots of unity, and by Petersson (\cite{P1} in characteristic $3$. In its full generality, however, Albert's question remains unresolved to this day. In the present paper, an affirmative answer will be given at the expense of passing to an isotope. More precisely, we establish the following result.
\vskip1mm
\noindent
    {\bf Theorem.} \emph{To every Albert division algebra over a field of arbitrary characteristic, there exists an isotope that contains a cyclic cubic extension of the base field}.
    \vskip1mm

    The proof, though surprisingly short, relies heavily on the theory of invariants, cohomological or otherwise, for cubic Jordan algebras. The paper concludes with a few immediate applications.

    This introduction would be incomplete without a mention of relations of the theory of Albert algebras with algebraic groups and, in particular, to the exceptional algebraic groups. We mention a few. 
    To begin with, for an Albert algebra $A$ defined over a field $k$ the full group $\text{\bf Aut}(A)$ of $k$-algebra automorphisms is a simple algebraic group of type $F_4$ defined over $k$ and all simple algebraic groups of type $F_4$ that are defined over $k$ arise in this fashion. In this, the split simple algebraic groups of type $F_4$ correspond to the split Albert algebras.

    Any Albert algebra $A$ carries a norm map, which is a cubic form $N$ defined on $A$. The full group $\text{\bf Str}(A)$ of similarities of $N$, called the structure group of $A$, is a reductive algebraic group of type $E_6$ defined over $k$ and its semisimple part is the full group $\text{\bf Isom}(A)$ of norm isometries of $A$, which is simple and simply connected of type $E_6$.

    A rank $2$ form of $E_8$ can be described in terms of Albert division algebras (see \cite{TW}), owing to the fact that such groups have anisotropic kernel the structure group of an Albert division algebra. Albert division algebras that are first Tits constructions contain cyclic cubic subfields (see \S {\bf 2.}{\bf c}-(iii)). This fact was exploited in (\cite{Th-1} and \cite{Th-2}) to prove the Tits-Weiss conjecture for the groups of type $E_8$ that arise from Albert division algebras that are first Tits constructions. We refer to (\cite{FF} and \cite{P0}) and references therein for further connections to Lie theory.

    Throughout this paper, we fix a base field $k$ of arbitrary characteristic. All Albert (and other) algebras are tacitly assumed to be over $k$. The reader may consult (\cite{KMRT}, Chap. IX and \cite{P0}) for their basic properties. For convenience, however, the ones that are most crucial for the subsequent applications will be summarized in the first few sections below. 

    \section{Preliminaries} In this section we briefly recall the definitions of the mod $2$ and mod $3$ Galois cohomological invariants of Albert algebras and some facts about the two rational construction of these algebras due to Tits. For details on Albert algebras and their cohomological invariants, we refer the reader to the survey article (\cite{P0}) or to (\cite{KMRT}). \\
    \vskip1mm
    \noindent
        {\bf (a).  The invariants mod $2$.} There are two types of Albert algebras: reduced ones and division algebras. A reduced Albert algebra $A$ may be realized by $3 \times 3$ twisted-hermitian matrices with entries in an octonion algebra $C$ and has two classifying invariants, the \emph{invariants mod $2$}, which are $n$-Pfister quadratic forms $f_n(A)$ for $n = 3,5$; in fact, $f_3(A)$ is just the norm of $C$. For $\text{char}(k) \neq 2$, $f_n(A)$ may be interpreted cohomologically as an element of $H^n(k,\mathbb{Z}/2\mathbb{Z})$ via the Arason invariant (see \cite{P0}, 7.5). Thanks to existence and uniqueness of the reduced model (\cite{P0}, 12.1), the invariants mod $2$ make sense also if $A$ is not reduced, though they are no longer classifying. \\
        \vskip1mm
        \noindent
        {\bf (b).  The invariant mod $3$.} For any Albert algebra $A$, there is another invariant, the \emph{invariant mod $3$}, denoted by $g_3(A)$. It belongs to $H^3(k,\mathbb{Z}/3\mathbb{Z})$ for $\text{char}(k) \neq 3$, and to the additive abelian group $H_3^3(k)$ defined in (\cite{PR4}, 2.) for $\text{char}(k) = 3$. It has the following fundamental properties.
\begin{itemize}
\item [(i)] $g_3(A)$ commutes with base change, i.e., $g_3(A \otimes l) = \text{res}_{l/k}(g_3(A))$, , where $\text{res}_{l/k}$ is the restriction homomorphism of Galois cohomology for $\text{char} \neq 3$, and its analogue explained in (\cite{PR4}, 3.) for $\text{char}(k) = 3$.

\item [(ii)] $g_3(A)$ detects division algebras: $A$ is an Albert division algebra if and only if $g_3(A) \neq 0$. 
\end{itemize}
It is not known whether the invariants mod $2$ and $3$ classify Albert algebras up to isomorphism.
\vskip1mm
\noindent
{\bf (c). The second Tits construction.}  The most versatile tool for constructing Albert algebras is the {second Tits construction}. Its \emph{input} is a quadruple $(B,\sigma,u,\mu)$, where $B$ is a central separable associative algebra over its centre $K$, a quadratic \'etale $k$-algebra, $\sigma :B \to B$ is a $K/k$-involution (so $(B,\sigma)$ is what we call a \emph{central simple associative $k$-algebra of degree $3$ with unitary involution}), and the pair $(u,\mu) \in B^\times \times K^\times$ is \emph{admissible} in the sense that $N_B(u) = N_K(\mu)$; here $N_B$ (resp. $N_K$) is the reduced norm of $B$ over $K$ (resp. the norm of $K$ as a quadratic $k$-algebra). The \emph{output} is an Albert algebra $A = J(B,\sigma,u,\mu)$ that contains
\[
A_0 := (B,\sigma)_+ := \{v \in B \mid \sigma(v) = v\}
\]
as a cubic Jordan subalgebra. The second Tits construction has the following fundamental properties.
\begin{itemize}
\item [(i)] $A$ is a division algebra if and only if $A_0$ is one and $\mu$ is not the reduced norm of an element in $B$.

\item [(ii)] By (\cite{PR5}, 3.5) and (\cite{PR4}, 8.), the invariant $\text{mod}~3$ of $A$ depends neither on $\sigma$ nor on $u$. 

\item [(iii)] If $K = k \oplus k$ splits, then, up to isomorphism, $B = D \oplus D^{\text{op}}$ for some central simple associative $k$-algebra $D$ of degree $3$ and $\sigma$ becomes the exchange involution. Moreover, if $u = 1_D \oplus 1_{D^{\text{op}}}$, then $\mu = \gamma \oplus \gamma^{-1}$ for some $\gamma \in k^\times$, and we abbreviate
\[
A = J(D,\gamma) := J(k \oplus k,D \oplus D^{\text{op}},1_D \oplus 1_{D^{\text{op}}},\gamma \oplus \gamma^{-1}),
\]
referring to this as the \emph{first Tits construction} arising from $D,\gamma$. It contains an isomorphic copy of $D^+$ as a cubic Jordan subalgebra.

\item [(iv)] Let $(B,\sigma)$ be a central simple associative $k$-algebra of degree $3$ with unitary involution. By (a), the Albert algebra $A := J(B,\sigma,1_B,1_K)$ is reduced. We put $f_3(\sigma) := f_3(A)$, dependence on $B$ being understood. The invariant $f_3$ classifies the unitary involutions on $B$ up to isomorphism (\cite{P2}, 2.4). We say $\sigma$ is \emph{distinguished} if $f_3(\sigma)$ is hyperbolic. If $A = J(B,\sigma,u,\mu)$ is any Albert algebra realized by the second Tits construction, then $f_3(A) = f_3(\sigma_u)$, where $\sigma_u$ is the unitary involution $v \mapsto u\sigma(v)u^{-1}$ on $B$ (\cite{PR6}, 1.8).
\end{itemize}
\section{\bf Proof of the theorem} In this section, we present a proof of the theorem on cyclicity of Albert division algebras and give some immediate applications.

\vskip1mm
\begin{proof} of the theorem. Let $A = J(B,\sigma,u,\mu)$ be any Albert division algebra over $k$ realized by means of the second Tits construction. By (\cite{KMRT}, (39.2)~(2)), whose proof obviously works in any characteristic, we may assume $N_B(u) = N_K(\mu) = 1$. Next we apply (\cite{P2}, 2.10) to exhibit an element $v \in B^\times$ such that $\sigma_v$ is distinguished.  Following
(\cite{Th}, 3.4), we now consider the Albert algebra $A^\prime := J(B,\sigma_v,1_B,\mu)$. From \textbf{2.}({\bf c})-(ii) we deduce $g_3(A^\prime) = g_3(A)$. On the other hand, $\sigma_v$ being distinguished and
  \textbf{2.}({\bf c})-(iv) imply that $f_3(A)$ is hyperbolic. By (\cite{P2}, 4.10), therefore, $A^\prime$ is a first Tits  construction. Combining Wedderburn's theorem on the cyclicity of associative division algebras of degree $3$ with \textbf{2.}({\bf c})-(iii), we conclude that $A^\prime$ contains a cyclic cubic subfield $L/k$. Hence $A^\prime_L$ is reduced (actually split), and
  \textbf{2.}({\bf b})-(i), (ii) imply $g_3(A_L) = 0$. Thus the Albert algebra $A_L$ is reduced, forcing $L$ to be a subfield of some isotope of $A$ (\cite{PR1}, Thm.~2).
  \end{proof}
\noindent
The following consequence was observed by Holger Petersson.
\vskip1mm
\noindent
{\bf Corollary.} \label{c.CYCRED} \emph{If $k$ has no cyclic cubic extensions, then every Albert algebra over $k$ is reduced.} \hfill $\square$ 

\noindent The preceding result generalizes (\cite{P4}, Thm.~7) and trivializes its proof. In particular, referring therein to the classification of Albert division algebras over local fields is no longer needed. 
\vskip1mm
\noindent
{\bf Corollary.} \label{c.CYCINV} \emph{Every Albert division algebra over $k$ has an isotope that contains a cyclic cubic subfield and whose $5$-invariant $\text{mod}~2$ is hyperbolic.}
\vskip1mm
\proof Let $A$ be an Albert division algebra over $k$. By the theorem, passing to an isotope if necessary, we may assume that $A$ contains a cyclic cubic subfield $L/k$. By (\cite{P2}, 4.7), there exists a $v \in L^\times$ such that the isotope $A^{(v)}$ has $f_5(A^{(v)})$ hyperbolic. Moreover, $L^{(v)} \cong L$ is a cyclic cubic subfield  of $A^{(v)}$. \hfill $\square$

\noindent The final result of the paper underscores the intimate connection between Albert algebras and exceptional groups.
\vskip1mm
\noindent
{\bf Corollary.} \label{CYCDEF} \emph{Let $A$ be an Albert division algebra over $k$ and assume $\text{char}(k) \neq 2, 3$. Then the structure group scheme of $A$ contains a subgroup of type $^3D_4$ defined over $k$.}
\vskip1mm
\proof The structure group scheme of $A$ does not change when passing to an isotope. By the theorem, therefore, we may assume that $A$ contains a cyclic cubic subfield. Applying (\cite{HT}, Cor.~3.2), we conclude that the automorphism group scheme of $A$ contains a subgroup of type $^3D_4$ defined over $k$. But the automorphism group scheme of $A$ is contained in its structure group scheme. \hfill $\square$

\vskip2mm
    \noindent
        {\bf Acknowledgement.} The author thanks Linus Kramer of Forschungs Institut Uni-Muenster for supporting this research in October 2018.
        The author thanks Holger P. Petersson with gratitude for his generosity and help in this work and Otmar Loos for many fruitful discussions during his visit to FernUni at Hagen.
        The suggestions made by the anonymous referee were invaluable and have improved the readability of the paper, we acknowledge his/her help.

\end{document}